\documentclass[12pt,a4paper]{article}

\usepackage[dvips]{graphicx}
\usepackage{calc}
\usepackage{color}
\usepackage{amsmath}
\usepackage{amssymb}
\usepackage{amscd}
\usepackage{amsthm}
\usepackage{amsbsy}
\usepackage{delarray}
\usepackage{enumerate}
\usepackage[T1]{fontenc}
\usepackage{inputenc}
\usepackage{enumerate}

\begin{document}



\setlength{\parindent}{5mm}
\renewcommand{\leq}{\leqslant}
\renewcommand{\geq}{\geqslant}
\newcommand{\N}{\mathbb{N}}
\newcommand{\sph}{\mathbb{S}}
\newcommand{\Z}{\mathbb{Z}}
\newcommand{\R}{\mathbb{R}}
\newcommand{\C}{\mathbb{C}}
\newcommand{\F}{\mathbb{F}}
\newcommand{\g}{\mathfrak{g}}
\newcommand{\h}{\mathfrak{h}}
\newcommand{\K}{\mathbb{K}}
\newcommand{\RN}{\mathbb{R}^{2n}}
\newcommand{\ci}{c^{\infty}}
\newcommand{\derive}[2]{\frac{\partial{#1}}{\partial{#2}}}
\renewcommand{\S}{\mathbb{S}}
\renewcommand{\H}{\mathbb{H}}
\newcommand{\eps}{\varepsilon}
\newcommand{\mumax}{\mu_{\max}}
\newcommand{\mumin}{\mu_{\min}}

\theoremstyle{plain}
\newtheorem{theo}{Theorem}[section]
\newtheorem{prop}[theo]{Proposition}
\newtheorem{lemma}[theo]{Lemma}
\newtheorem{definition}[theo]{Definition}
\newtheorem*{notation*}{Notation}
\newtheorem*{notations*}{Notations}
\newtheorem{corol}[theo]{Corollary}
\newtheorem{conj}[theo]{Conjecture}
\newtheorem{question}[theo]{Question}
\newtheorem*{question*}{Question}
\newenvironment{demo}[1][]{\addvspace{8mm} \emph{Proof #1.
    ---~~}}{~~~$\Box$\bigskip}

\newlength{\espaceavantspecialthm}
\newlength{\espaceapresspecialthm}
\setlength{\espaceavantspecialthm}{\topsep-3pt} \setlength{\espaceapresspecialthm}{\topsep-13pt}

\newenvironment{example}[1][]{\refstepcounter{theo} 
\vskip \espaceavantspecialthm \noindent \textsc{Example~\thetheo
#1.} }%
{\vskip \espaceapresspecialthm}

\newenvironment{remark}[1][]{\refstepcounter{theo} 
\vskip \espaceavantspecialthm \noindent \textsc{Remark~\thetheo
#1.} }%
{\vskip \espaceapresspecialthm}

\def\Homeo{\mathrm{Homeo}}
\def\Hameo{\mathrm{Hameo}}
\def\Diffeo{\mathrm{Diffeo}}
\def\Symp{\mathrm{Symp}}
\def\Id{\mathrm{Id}}
\newcommand{\norm}[1]{||#1||}
\def\Ham{\mathrm{Ham}}
\def\Cal{\mathrm{Cal}}

\title{Hofer's distance on diameters and the Maslov index}
\author{Vincent Humili\`ere$^1$}
\normalsize  \footnotetext[1]{Institut Mathématique de Jussieu, Université Pierre et Marie Curie, Paris, France

Email: \texttt{humi@math.jussieu.fr}}
\maketitle

\abstract{We prove that Hofer's distance between two diameters of the open 2-disk admits an upper bound in terms of the Maslov index of their intersection points.}

\section{Introduction and result}

This note is inspired by Khanevsky's article \cite{khanevsky}. We will adopt his settings throughout the paper. Let us recall them briefly. 
Let $D\subset\R^2=\{(x,y)\,|\,x\in\R,y\in\R\}$ be the open unit 2-disk, endowed with the symplectic structure $\omega=\frac1\pi dx\wedge dy$ (so that the total area of $D$ is 1). We call \emph{standard diameter} and denote by $L_0$ the intersection of $D$ with the $x$-axis $\R\times\{0\}$. Given a time-dependent Hamiltonian $H$ on $D$, we denote by $X_H$ its Hamiltonian vector field defined by $\omega(X_H,\cdot)=dH$, and by $\phi_H^t$ the Hamiltonian isotopy generated by $H$. All Hamiltonian functions will be assumed compactly supported. By definition, a \emph{diameter} is a curve which is isotopic to the standard diameter via a compactly supported Hamiltonian isotopy. We denote by $\mathcal{E}$ the set of all diameters. 

The set $\mathcal{E}$ can be endowed with the so-called \emph{Hofer distance} which is defined as follows:
$$d(L_1,L_2)=\inf\left\{\int_0^1\left(\max_{x\in D} H_t-\min_{x\in D} H_t\right)dt\right\},$$
where the infimum is taken over all smooth and compactly supported time dependent Hamiltonians $H_t$ such that $\phi_H^1(L_1)=L_2$. The distance $d$ may be defined on more general symplectic manifolds. Its non-degeneracy is a deep result (see \cite{chekanov}). The value of  $\int_0^1(\max {H_t}-\min {H_t})dt$ is called the \emph{Hofer energy} of the isotopy generated by $H_t$.

An interesting fact pointed out in Khanevsky's paper is that Hofer's distance descends to a metric on the \emph{reduced diameter space} which is defined as 
$\hat{\mathcal{E}}= \mathcal{E}/\sim$, where $L\sim L'$ iff $L=\phi(L')$ with $\phi\in\mathcal{S}$, $\mathcal{S}$ being the group of Hamiltonian diffeomorphisms which globally fix the standard diameter $L_0$. Several other invariants of $\mathcal{E}$ descend to $\hat{\mathcal{E}}$. Khanevsky considers two of them: the number of intersection points between diameters and (much more involved) an invariant $r_A$ constructed from the Entov-Polterovich quasimorphism and the Calabi morphism. Roughly speaking, he proves that $r_A$ is dominated by $d$ (this implies that $\hat{\mathcal{E}}$ is unbounded) and that $d$ is dominated by the number of intersection points. In the present note we are interested in a third invariant which is the Maslov index.

\medskip
The Maslov index of a transverse intersection point of two given Lagrangians is an integer which is invariant under Hamiltonian transformations. Unless we choose some convention, this index is defined up to an additive constant but the difference of index between two points (the "index gap") is well defined. It can be constructed in a natural way by an abstract construction \cite{viterbo}. However, in our case, it can be constructed in a naive way as follows.

First, applying a Hamiltonian diffeomorphism, we may suppose that one of our two diameters is the standard one. Then, consider a diameter $L$ intersecting $L_0$ transversely. Note that in our context, every diameter coincides with the standard diameter in the complement of a compact set. Therefore, we say that $L$ and $L_0$ are transverse if we can write $L\cap L_0=(-1,-1+\eps]\cup\{x_i\}_{i\in \{1,\ldots,N\}}\cup[1-\eps',1)$, where $\eps,\eps'>0$ and the points $x_1,\ldots,x_N$ are transverse intersections (in the usual sense). We suppose that the points $x_1,\ldots,x_N$ are ordered by their position on $L$ (not $L_0$). We choose a convention for the index of $x_1$, namely we set $\mu(x_1)=0$. Then, we construct the Maslov index $\mu$ of the other intersection points inductively as follows. 

We denote by $D^+$ the upper half-disk and $D^-$ the lower half-disk. For $i\in\{1,\cdots,N\}$, we set $\mu(x_{i+1})=\mu(x_i)+\delta$, where
$$\delta=\left\lbrace
\begin{array}{cl}
  1, & \text{if }L|_{[x_i,x_{i+1}]}\subset D^-\text{ and }x_i < x_{i+1}\text{ on } L_0\\
  -1, & \text{if }L|_{[x_i,x_{i+1}]}\subset D^-\text{ and }x_i > x_{i+1}\text{ on } L_0\\
  -1, & \text{if }L|_{[x_i,x_{i+1}]}\subset D^+\text{ and }x_i < x_{i+1}\text{ on } L_0\\
  1, & \text{if }L|_{[x_i,x_{i+1}]}\subset D^+\text{ and }x_i > x_{i+1}\text{ on } L_0                     
\end{array}
\right.$$

\begin{remark}\label{remarque indice}\begin{enumerate}
\item Intuitively, $\mu(x_i)$ measures how much $L$ "twists" in the positive direction before reaching $x_i$.
\item The index $\mu$ descends to the reduced diameter space $\hat{\mathcal{E}}$. 
\item Note that by definition, the index gap between two intersection points that are consecutive on $L$ is $1$ or $-1$. The abstract definition of the Maslov index implies that we can invert the roles of $L$ and $L_0$ and therefore that the index gap between two intersection points that are consecutive on $L_0$ is also either $1$ or $-1$.
\end{enumerate}\end{remark}

To avoid any confusion, we will sometimes denote $\mu(x,L)$ instead of $\mu(x)$. We then set
$$\mu_{\max}(L)=\max_{i=1,\ldots,N}\mu(x_i,L)\quad\text{ and }\quad\mu_{\min}(L)=\min_{i=1,\ldots,N}\mu(x_i,L).$$
The maximal index gap $\mu_{\max}(L)-\mu_{\min}(L)$ descends to the reduced diameter space $\hat{\mathcal{E}}$. It is therefore natural to try to compare it with the other invariants considered by Khanevsky. The following inequality is a consequence of Remark \ref{remarque indice}.4.
$$2(\mu_{\max}(L)-\mu_{\min}(L))\leq \sharp(L\cap L_0)+1.$$
Our main result is then the following theorem.

\begin{theo}\label{theoreme hofer-maslov} For any diameter $L$ transverse to $L_0$, with at least two transverse intersection points with $L_0$,
$$d(L,L_0)\leq \mumax(L)-\mumin(L)-\frac12.$$
\end{theo}

\begin{remark} The maximal index gap $\mumax(L)-\mumin(L)$ vanishes if and only if there is only one transverse intersection point. In that case, $d(L,L_0)\leq \frac12$.
\end{remark}

\begin{remark} This implies $d(L,L_0)\leq \frac12\cdot\sharp(L\cap L_0)$. Therefore our result implies a linear estimate as in \cite{khanevsky} but with $\frac12$ as multiplicative constant, while Khanevsky's constant is $\frac18$. Nevertheless, in many cases, our upper bound is much better. Moreover, after I completed this note, Michael Khanevsky explained to me an argument based on his paper's settings (in particular "diameter trees") showing that it might be possible to improve the multiplicative constant 1 in Theorem \ref{theoreme hofer-maslov} to the constant $\frac14$.\end{remark}

\begin{remark} As far as I know, the analogous problem in higher dimension is completely open.
\end{remark}

Let us give an idea of the proof of Theorem \ref{theoreme hofer-maslov}. Suppose that $\mumin(L)<0$. We first concentrate on intersection points of minimal index. Given such a point, it is possible to construct a very simple Hamiltonian isotopy with several properties including the fact that it removes the intersection point (see Section 3). We then show that such simple transformations can be performed in an appropriate order so that all intersection points of minimal index are removed and so that the total energy needed for these transformation is less than 1. We get a new diameter $L'$ satisfying $\mumin(L')\leq \mumin(L)+1$ satisfying $d(L',L)=1$. We then proceed by induction until $\mumin=0$, and work similarly for points of positive index (see Section 4).

\begin{remark} A consequence of our proof is that for any diameter $L$ with $\mumin(L)=0$, there exists a non-negative Hamiltonian function $H$ such that $d(\phi_H^1(L_0),L)\leq 1$. 
\end{remark}

\subsection*{Aknowlegments} I am very grateful to Frédéric Le Roux for useful discussions and for listenning to me patiently. I also thank Leonid Polterovich for drawing my attention on Michael Khanevsky's paper. I thank the three of them for their comments and suggestions on the first version of the paper. I am also grateful to H\'el\`ene Eynard for her careful reading. Finally I thank an anonymous referee for pointing out some mistakes and for his suggestions to improve the paper.

\section{Two lemmas on the Maslov index}

Let $L$ be a fixed diameter. We first recall the abstract definition of the Maslov index gap between two transverse intersection points, as introduced by Viterbo in \cite{viterbo}. Let $x,y$ be two transverse intersection points. Denote by $\gamma:[0,1]\to L$ a curve on $L$ with $\gamma(0)=x$ and $\gamma(1)=y$. The derivative $\dot{\gamma}$ determines a path $[\dot{\gamma}]$ in the projective space $\R P^1=\R/\Z$. This path can be turned into a loop after concatenation with a path $\tau:[0,1]\to\R P^1$, such that $\tau(0)=[\dot{\gamma}](1)$, $\tau(1)=[\dot{\gamma}](0)$ and such that for any $s$, $\tau(s)$ is transverse to $L_0$.  Let $\tilde{\gamma}:[0,1]\to\R$ be the lift of this loop satisfying $\tilde{\gamma}(0)=0$. The index gap $\mu(y,L)-\mu(x,L)$ is then the Maslov index of the loop, i.e, the integer $\tilde{\gamma}(1)$.

We will need the two following lemmas.

\begin{lemma}\label{lemme indices preserves} Let $x,y$ be two transverse intersection points between $L$ and $L_0$. Let $(h^t)_{t\in[0,1]}$ be a compactly supported Hamiltonian isotopy of the disk. Suppose that for any $t\in[0,1]$, the points $x$ and $y$ remain intersection points of $L_0$ and $h^t(L)$. Suppose moreover that the tangent line to $h^t(L)$ at $x$ (resp. $y$) remains constant for $t\in[0,1]$. Then the index gap remains constant:
$$\forall t\in[0,1],\  \mu(y,h^t(L))-\mu(x,h^t(L))=\mu(x,L)-\mu(y,L).$$
\end{lemma}
\begin{demo} Let $\gamma$ be a curve joining $x$ to $y$ on $L$, as above. Let $\alpha_t$ be a continuous family of curves joining $x$ to $h^t(x)$ on $h^t(L)$ and $\beta_t$ a continuous family of curves joining $h^t(y)$ to $y$ on $h^t(L)$. Then consider the family of curves $\gamma_t$ obtained as a smooth reparametrization of the concatenation of three paths:
$$\gamma_t=\alpha_t\ast (h^t\circ\gamma)\ast\beta_t.$$
Our assumptions imply that the family of loops $[\dot{\gamma_t}]\ast\tau$, where $\tau$ is as above, is a homotopy with fixed end points in $\R P^1$. Thus,  $\tilde{\gamma_t}(1)=\tilde{\gamma}(1)$ for any $t$, and the lemma follows. 
\end{demo}

As in the introduction, we write $L\cap L_0=(-1,-1+\eps]\cup\{x_i\}_{i\in \{1,\ldots,N\}}\cup[1-\eps',1)$, where $\eps,\eps'>0$ and the points $x_1,\ldots,x_N$ are transverse intersections ordered by their position on $L$.

\begin{lemma}\label{lemme x1xN} We suppose that both arcs $L|_{[0,x_1]}$ and $L|_{[x_N,1]}$ are included in $D^+$. Then, 
$\mu(x_N)-\mu(x_1)=1$.
\end{lemma}
\begin{figure}
\begin{center}
\includegraphics[width=7cm]{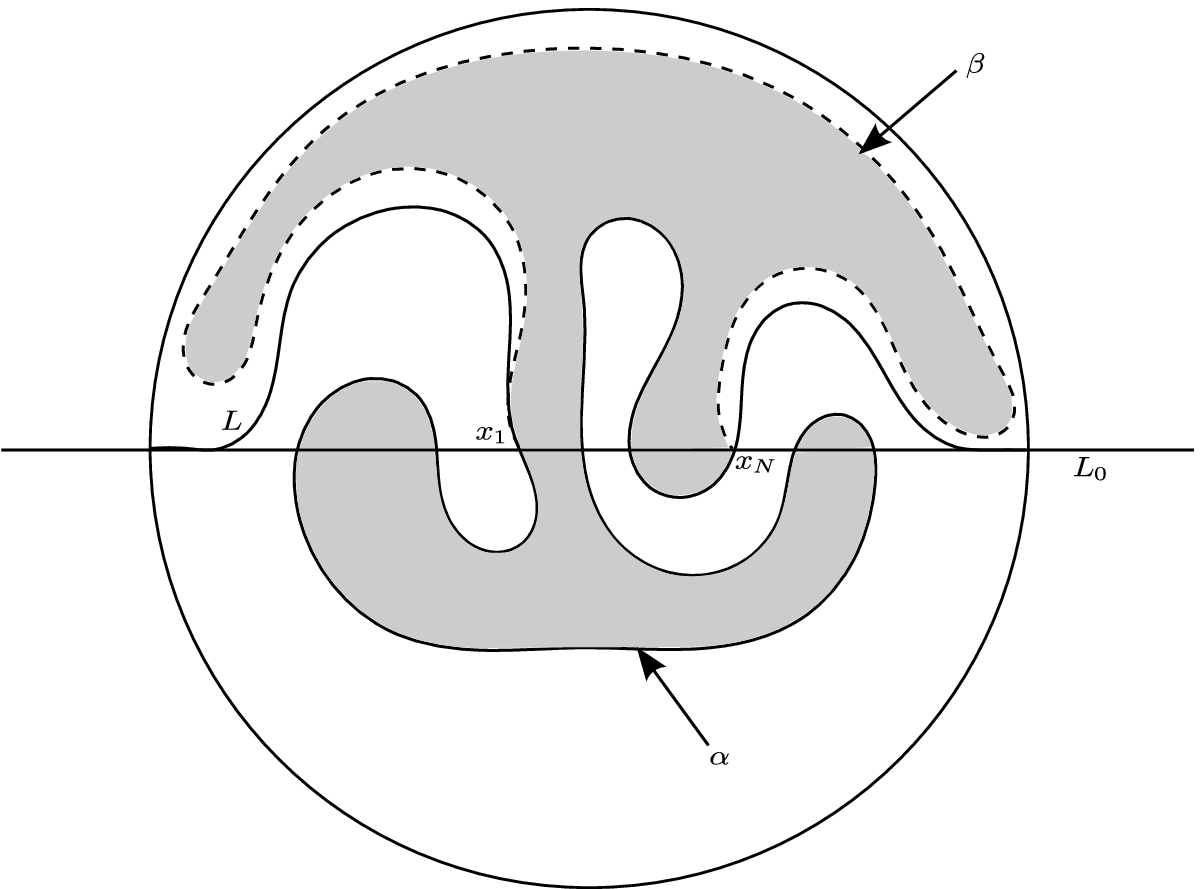} 
\end{center}
\caption{Proof of Lemma \ref{lemme x1xN}}
\label{dessin x1xN} 
\end{figure} 
\begin{demo} Let $\psi$ be a Hamiltonian diffeomorphism such that $L=\psi(L_0)$. Let $\gamma$ be a smooth path on $D$ which is the concatenation of two paths: a path $\alpha$ that parametrizes the arc $L_{[x_1,x_N]}$ and a path $\beta$ from $x_N$ to $x_1$ disjoint from $L_0$ and $L$ (see Figure \ref{dessin x1xN}). The path $\gamma$ is a simple curve which is the boundary of an embeded disk. Therefore, the Maslov index of $[\dot{\gamma}]$ equals 2. Since the contribution of $\beta$ to the Maslov index is obviously 1, it follows that the contribution of $\alpha$ is also 1. Thus, $\mu(x_N)-\mu(x_1)=1$. 
\end{demo}

\section{Removing points of minimal index}

Let $L$ be a diameter transverse to $L_0$. In this section we describe simple Hamiltonian transformations that allow to remove an intersection point with extremal index. We describe it only for points of minimal index. Points of maximal index can be treated in a similar way.

\begin{figure}
\begin{center}
\includegraphics[width=6cm]{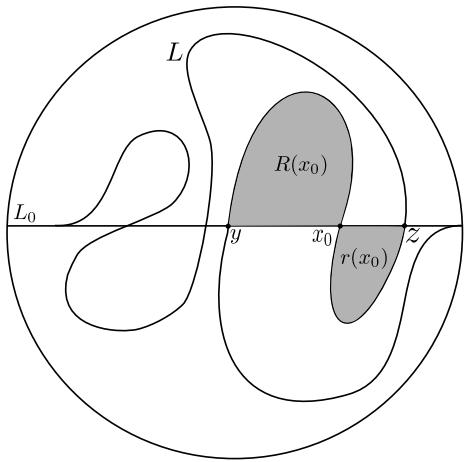} 
\end{center}
\caption{Large and small regions of a minimal index point}
\label{dessin petites_grandes_regions} 
\end{figure} 

We start with the following remark on the shape of $L$ near a point of minimal index. Let $x_0$ be a point with minimal Maslov index. Consider the two  transverse intersection points of $L$ and $L_0$ which are consecutive to $x_0$ on $L$. Here, we make the assumption that those two points exist, i.e., that $x_0\neq x_1$ and $x_0\neq x_N$. The minimality of $\mu(x)$ imposes that one of them (call it $y$) is on the left of $x_0$ on $L_0$ and the other (call it $z$) is on its right. It also forces the arc $\alpha=L|_{[y,x_0]}$ to lie in $D^+$ and the arc $\beta=L|_{[x_0,z]}$ to lie in $D^-$ (see figure \ref{dessin petites_grandes_regions}).

We assume that one of the two closed regions delimited respectively by $\alpha$ and $L_0$ on one side, $\beta$ and $L_0$ on the other side has a smaller area. This region will be called the \emph{small region} of $x_0$ and denoted by $r(x_0)$. Its area will be called the \emph{weight} of $x_0$ and denoted by $w(x_0)$. The largest region will be called the \emph{large region} of $x_0$ and denoted $R(x_0)$. By \emph{half-disk} we will mean any intersection of a smooth closed disk with one of the two standard half disks $D^+$ and $D^-$ (figure \ref{dessin petites_grandes_regions}).

\begin{figure}
\begin{center}
\includegraphics[width=13cm]{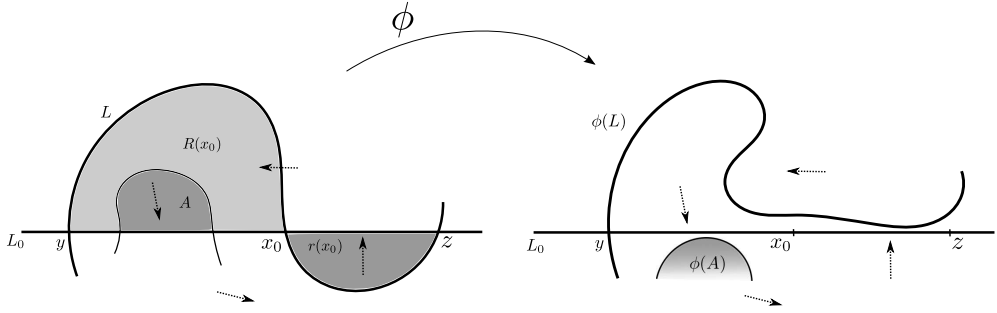} 
\end{center}
\caption{Illustration of Lemma \ref{lemme technique} (the dotted arrows indicate roughly how $\phi$ acts)}
\label{dessin lemme technique} 
\end{figure}

\begin{lemma}[Removing a minimal index point, figure \ref{dessin lemme technique}]\label{lemme technique}
Let $\eps>0$, $U$ a neighbourhood of $r(x_0)\cup R(x_0)$, and $A$ any half disk included in the interior of $R(x_0)$, of area less than $w(x_0)$ and such that any connected component of $A\cap L$ is an arc with at least one extremity on $L_0$. Then there exists a 
Hamiltonian diffeomorphism $\phi$ supported in $U$ with $\|\phi\|\leq w(x_0)+\eps$, such that
 \begin{itemize}
  \item $\phi$ removes at least the intersection points $x_0$ and $z$ and does not create new intersection points, i.e., $\phi(L)\cap L_0\subset L\cap L_0-\{x_0,z\}$,
  \item $\phi$ does not change the index gap between any two remaining intersection points,
  \item $\phi$ maps the small region $r(x_0)$ into the opposite half disk. In particular it removes all intersection points contained in $r(x_0)$.
  \item $\phi$ maps $A$ into the opposite half disk.
 \end{itemize}
\end{lemma}

\begin{figure}
\begin{center}
\includegraphics[width=6cm]{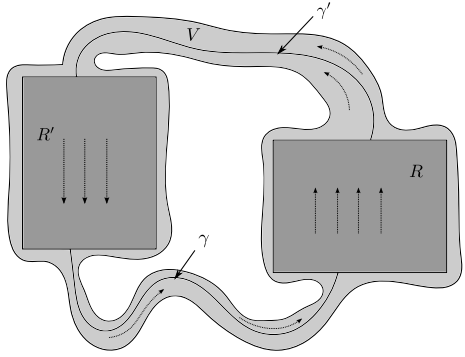} 
\end{center}
\caption{Illustration of Lemma \ref{lemme hamiltonien precis}.}
\label{dessin hamiltonien simple} 
\end{figure}

We will use several times the following fact and therefore state it in a separate lemma (see Figure \ref{dessin hamiltonien simple}).
 
\begin{lemma}\label{lemme hamiltonien precis} Let $R=[a,b]\times[c,d]$, $R'=[a',b']\times[c',d']$ be two disjoint closed rectangles in $\R^2$ with same area $\alpha$. We suppose moreover that $a'<b'<a<b$. Let $\gamma$ (resp. $\gamma'$) be a curve that does not meet $R\cup R'$ and joins a point of $[a,b]\times \{c\}$ (resp. $[a,b]\times \{d\}$) to a point of $[a',b']\times \{c'\}$ (resp. $[a',b']\times \{d'\}$). We also suppose that both curves are disjoint.
Let $\eps>0$ and let $V$ be an open neighbourhood of $R\cup R'\cup\gamma\cup\gamma'$. Then there exists a smooth autonomous Hamiltonian function $H$ satisfying
\begin{itemize}
 \item $H=0$ on the unbounded component of $\R^2-V$ and $H=\alpha+2\eps$ on the bounded component of $\R^2-V$.
 \item $\min H =0$ and $\max H=\alpha+2\eps$
 \item On each rectangle $R$ and $R'$, $H$ is an affine function of $x$. More precisely, $H(x,y)=(c-d+\eps)(x-b)+\eps$ on $R$ and $H(x,y)=(d'-c'+\eps)(x-a')+\eps$ on $R'$.
 \item $\phi_H^1(R)=R'$ and $\phi_H^1(R')=R$.  
\end{itemize}
\end{lemma}

Note that the second condition on $H$ implies that the Hofer energy of the time one map $\phi_H^1$ is less or equal to $\alpha+2\eps$. Note also that the third condition means that the flow lines restricted to both rectangles are vertical straight lines.

\begin{figure}
\begin{center}
\includegraphics[width=6cm]{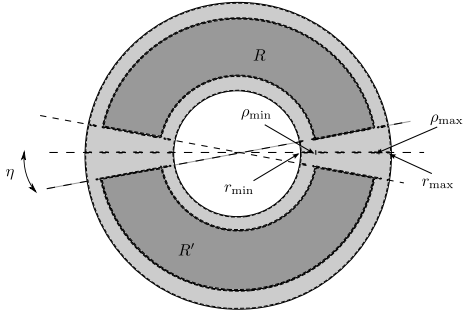} 
\end{center}
\caption{Image of Figure \ref{dessin hamiltonien simple} by some area preserving diffeomorphism.}
\label{dessin anneau} 
\end{figure} 

\begin{proof} Under some area preserving map, the situation is conjugated to that where $V$ is a standard annulus of the form $\{(x,y)\in\R^2\,|\,r_{\min}^2<x^2+y^2<r_{\max}^2\}$. The transformation can also be chosen so that the rectangles become angular sectors of a sub-annulus. In polar coordinates $(\rho,\theta)$, we may suppose that $R$ is now $$\{(\rho,\theta)\in(0,+\infty)\times(-\pi,\pi]\,|\,\rho_{\min}\leq \rho \leq\rho_{\max}\text{ and } \eta\leq\theta\leq \pi-\eta\}$$ and that $R'$ is $$\{(\rho,\theta)\in(0,+\infty)\times(-\pi,\pi]\,|\,\rho_{\min}\leq \rho \leq\rho_{\max}\text{ and } -\pi+\eta\leq\theta\leq -\eta\},$$
for some real numbers $\rho_{\min}$, $\rho_{\max}$ and $\eta$ such that $r_{\min}<\rho_{\min}<\rho_{\max}<r_{\max}$, and $0<\eta$.
The straight vertical lines in the rectangles become pieces of circles centered at the origin. See figure \ref{dessin anneau}.

Therefore, we only need to find $H$ in this situation, which is easy: take $H$ in the form $H(\rho,\theta)=h(\rho)$, where $h$ is some non-increasing function which satisfies the following properties.
\begin{itemize}
 \item if $\rho\leq r_{\min}$ then $h(\rho)=\alpha+2\eps$,
 \item if $\rho_{\min}\geq\rho\geq \rho_{\max}$ then $h(\rho)=\frac12(\rho_{\max}^2-\rho^2)+\eps$,
 \item if $\rho\geq r_{\max}$ then $h(\rho)=0$.
\end{itemize}
The associated flow acts as a linear rotation on the sub-annulus delimited by $\rho_{\min}$ and $\rho_{\max}$. The time one map makes half a turn and hence sends $R$ onto $R'$ and $R'$ onto $R$ as wished.
\end{proof}

\begin{demo}[of Lemma \ref{lemme technique}] First, recall that at any moment of the proof we can apply a small Hamiltonian perturbation and that all objects may thus be supposed to be in generic position. Let $\eps>0$ sufficiently small and $U$ an open neighbourhood of $r(x_0)\cup R(x_0)$. The first step of the proof is to note that after conjugation by some area preserving map of the plane, we may suppose that $D^- $ and $D^+$ are rectangles, that $r(x_0)$ is included in a rectangle $R$ af area less than $w(x_0)+\eps$ and that this rectangle $R$ is included in $U\cap D^-$. We may also suppose that $A$ lies inside a rectangle $R'\subset D^+$ included in $R(x_0)$ and with the same area as $R$, with $R\cap R'=\emptyset$. The arcs which are the connected components of $L\cap R'$ either have both extremities on $L_0$ or may be assumed to be pieces of straight vertical lines with one extremity on $L_0$.

In this situation, let $\gamma'$ be a smooth curve included in $D^+\cap U$ joining (as in Lemma \ref{lemme hamiltonien precis}) the top sides of $R$ and $R'$, and $\gamma$ be  a smooth curve in $D^-\cap U$ joining the bottom sides of $R$ and $R'$. Let $V$ be a neighborhood of $R\cup R'\cup\gamma\cup\gamma'$ included in $U$ and sufficiently narrow to avoid meeting $L_0$ except in the neighbourhood of $R'$ and $R$. Now, we may apply Lemma \ref{lemme hamiltonien precis} which provides us with a Hamiltonian function $H$. Then set $\phi=\phi_H^{\tau}$ where $\tau<1$ will be appropriately chosen in the sequel.

Since $\phi_H^1$ sends $R$ onto $R'$, $r(x_0)$ is sent entirely into $D^+$ by $\phi$ provided $1>\tau>\inf\{t\,|\,\phi_H^t(r(x_0))\subset D^+\}$. In particular the intersection points $x_0$ and $z$ disappear. In the same way $A$ is sent into $D^-$. We have to make sure that no new intersection point appears. These new points would necessarily come from the arcs of $L$ meeting $V$. If $V$ is sufficiently narrow, such an arc must meet $R$, $R'$ or $\gamma\cup\gamma'$. First remark that such an arc cannot meet $R$ except if it is included in $r(x_0)$ ($R$ is supposed sufficiently narrow). But then, the whole arc is sent to $D^+$ and no new intersection point appears. Then, note that the assumption made on the intersection $R'\cap L$ implies that along the isotopy, the intersection points either disappear or do not move at all and no new point appears. Finally, the arcs which correspond to intersections between $L$ and $\gamma\cup\gamma'$ do not create new intersection points either provided $\tau$ is sufficiently close to $\inf\{t\,|\,\phi_H^t(r(x_0))\subset D^+\}$.

The intersection points that remain after applying $\phi$ are either those not contained in the support of $\phi$ or some of those included in $R'$. In both cases these points satisfy the assumptions of Lemma \ref{lemme indices preserves}. Thus, for any pair of remaining points $u$, $v$, Lemma \ref{lemme indices preserves} gives: $\mu(v,\phi(L))-\mu(u,\phi(L))=\mu(v,L)-\mu(u,L)$.

Finally, according to Lemma \ref{lemme hamiltonien precis}, our diffeomorphism $\phi$ has energy less than approximately the area of the rectangle $R$. But this area is also approximately the area of $r(x_0)$, that is $w(x_0)$.  
\end{demo}

\begin{figure}
\begin{center}
\includegraphics[width=9cm]{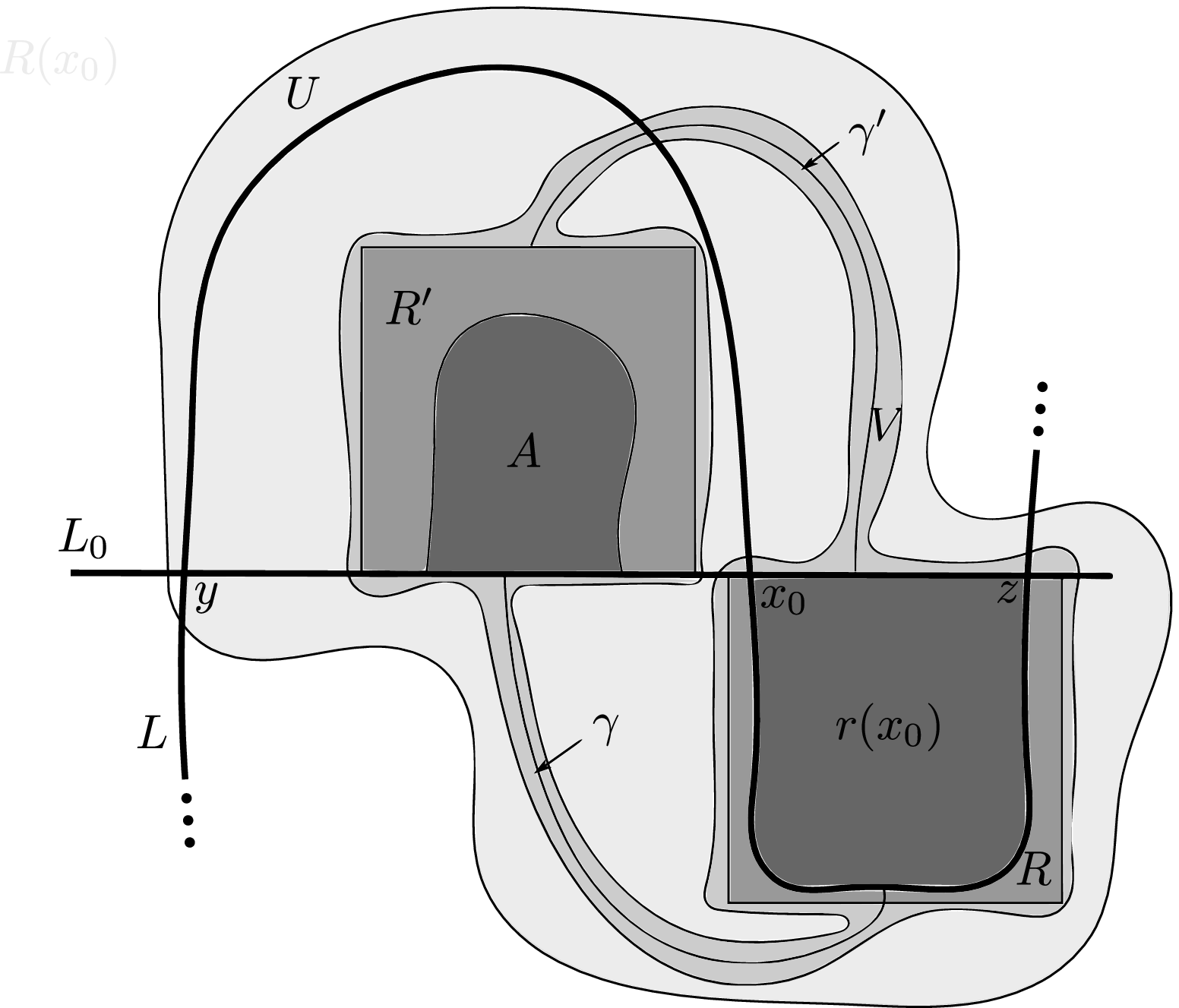} 
\end{center}
\caption{Proof of Lemma \ref{lemme technique}.}
\label{dessin preuve lemme technique} 
\end{figure} 


\begin{remark}\label{remark partial removing} Note that if we take the isotopy $\phi_H^t$ of the proof of Lemma \ref{lemme technique} at some smaller time, we can construct for any half-disk $A\subset R(x_0)$ of area less than $w(x_0)$ a diffeomorphism similar to that in Lemma \ref{lemme technique} that does not remove the intersection point $x_0$, but decreases $r(x_0)$ by approximately the area of $A$. The energy needed is approximately the area of $A$.
We will call this procedure a \emph{partial removing} of $x_0$ associated to $A$.
\end{remark}

\begin{figure}
\begin{center}
\includegraphics[width=10cm]{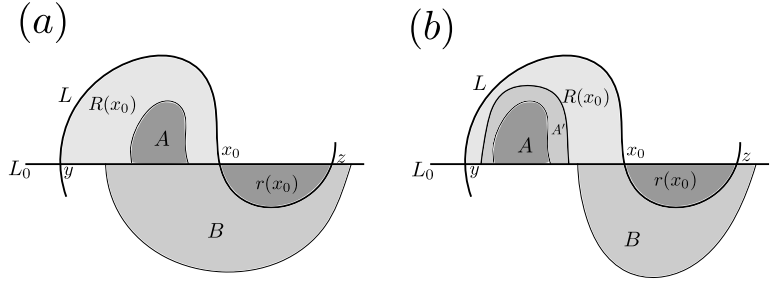} 
\end{center}
\caption{Illustration of Remark \ref{remark region du bas}}
\label{dessin remarques} 
\end{figure}

\begin{remark}\label{remark region du bas} We consider again an isotopy like the one appearing in the proof of Lemma \ref{lemme technique}. Let $B$ be a half disk included in $D^-$ and containing $r(x_0)$. It will be useful to figure out how the area of $B$ behaves along the isotopy $\phi_H^t$ in the following two cases. 

First, if $A\cap L_0\subset B\cap L_0$, then it is possible to choose the support of the isotopy so that it does not meet $\partial B\backslash L_0$. In this case, the area of $B$ remains constant along the isotopy (see Figure \ref{dessin remarques} (a)). 

Second, if $A$ is entirely on the left of $B$ (i.e., $\forall a\in A\cap L_0,\forall b\in B\cap L_0, a\leq b$). Then, for any other half disk $A'$ containing $A$ and on the left of $B$, we can choose the path $\gamma$ (resp. $\gamma'$) so that once it has entered $B$ (resp. $A'$) it does not exit $B$ (resp. $A'$) anymore (i.e., $\gamma\cap B$ is connected). This choice made, we see that along our (area preserving) isotopy, the respective areas of $A$, $A'$ and $B$ decrease at the same speed exactly (see Figure \ref{dessin remarques} (b)).
\end{remark}

\begin{figure}
\begin{center}
\includegraphics[width=9cm]{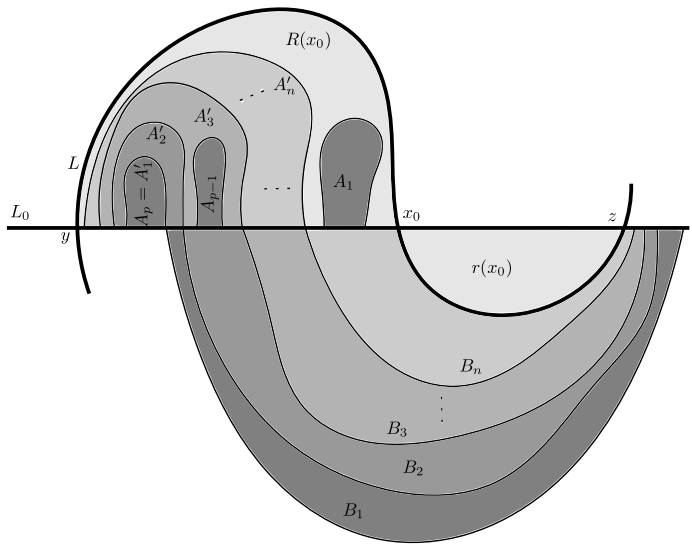} 
\end{center}
\caption{Settings of Lemma \ref{lemme technique 2}}
\label{dessin lemme technique 2} 
\end{figure}

With the help of these remarks, we can prove the following refinement of Lemma \ref{lemme technique}, illustrated by Figure \ref{dessin lemme technique 2}.

 \begin{lemma}\label{lemme technique 2} Let $\eps>0$ and $U$ be an open neighbourhood of $r(x_0)\cup R(x_0)$.
Let $A_1,\ldots, A_p=A'_1\ldots A'_n, B_1,\ldots, B_n$ be closed half disks such that
\begin{itemize}
 \item the $A_i$'s are included in $R(x_0)$, are pairwise disjoint and are ordered by their position on $L_0$ from the right to the left,
 \item $\sum_{i=1}^{p-1}\textrm{area}(A_i)<w(x_0)$,
 \item $A_p=A'_1\subset \cdots\subset A'_n\subset R(x_0)\subset D_+$,
 \item $r(x_0)\subset B_n\subset\cdots\subset B_1=B\subset D_-$,
 \item for any index $i$, $A_i'$ is on the left of $B_i$ and intersects at one point exactly,
 \item for any indices $i,j$, either $A_i\subset A_j'$ or $A_i\cap A_j'=\emptyset$.
\end{itemize}
Then, there exists a 
Hamiltonian diffeomorphism $\phi$ supported in $U$ with $\|\phi\|\leq w(x_0)+\eps$, such that
 \begin{itemize}
  \item $\phi(L)\cap L_0\subset L\cap L_0-\{x_0,z\}$,
  \item $\phi$ does not change the index gap between two remaining intersection points,
  \item $\phi$ maps the small region $r(x_0)$ into the opposite half disk, 
  \item $\phi$ maps each $A_i$, $1\leq i\leq p-1$ into the opposite half disk,
  \item for any index $j\in\{1,\ldots, n\}$, the map $\phi$ decreases the areas of $A'_j$ and $B_j$ by exactly the same amount, except in the special case $j=1$ and  $w_0>\sum_{i=1}^{p}\textrm{area}(A_i)$, for which $A_1'=A_p$ is entirely mapped to the opposite half disk.
 \end{itemize}
\end{lemma}

\begin{demo} We first apply successively partial removings of $x_0$ associated to each $A_i$ with $1\leq i\leq p-1$ (see Remark \ref{remark partial removing}). We get diffeomorphisms $\phi_1$,...,$\phi_{p-1}$.  We can choose them so that the support of each $\phi_i$ meets neither the $A_j$'s for $j>i$, nor the $A'_j$'s that do not contain $A_i$. Moreover, according to Remark \ref{remark region du bas}, this can be done so that for any index $j$, if $A'_j$ contains $A_i$, then the areas of $A'_j$ and $B_j$ are reduced by the same amount exactly, and if $A'_j$ does not contain $A_i$, then $A'_j$ and $B_j$ remain globally unchanged.

After the action of the composition $\phi_{p-1}\circ\cdots\circ\phi_{1}$, the half-disks $A_1,\cdots,A_{p-1}$ have been moved to $D_-$, the small region $r(x_0)$ has now area $c=w(x_0)-\sum_{i=1}^{p-1}\textrm{area}(A_i)$ and moreover $A_p=A'_1$ and $B_1$ remain unchanged. 

Then, two cases may occur. If $c>\textrm{area}(A_p)$, then we apply Lemma \ref{lemme technique} for a half-disk $A$ containing $A_p$ and of area approximately $c$. If $c<\textrm{area}(A_p)$, then we apply Lemma \ref{lemme technique} for a half-disk $A$ included in $A_p$ and of area approximately $c$. In the first case, this gives a diffeomorphism $\phi_p$ which maps $A_p$ to $D_-$. In the second case, this gives a diffeomorphism $\phi_p$ which, according to Remark \ref{remark region du bas}, can be chosen so that both areas of $A_p=A'_1$ and $B_1$ are decreased by $c$. Finally we see that the diffeomorphism  $\phi=\phi_{p}\circ\cdots\circ\phi_{1}$ suits our needs.
\end{demo}

\section{Proof of the theorem}

We will prove Theorem \ref{theoreme hofer-maslov} by induction on the maximal index gap $\mumax(L)-\mumin(L)$. Before to start, we apply an arbitrarily small perturbation on $L$ so that we may assume that $\mu(x_1)=\mu(x_N)=0$. This can be achieved by creating new transverse intersection points near the boundary to be in the situation of Lemma \ref{lemme x1xN} and then add one more intersection point to get $\mu(x_1)=\mu(x_N)$.

First, suppose that $\mumax(L)-\mumin(L)=1$. In this case, the intersection points of $L$ and $L_0$ are in a very simple configuration: the indices are successively $0$, $\pm 1$, 0, $\pm1$... and there exists some area preserving map of the plane, that preserves the $x$-axis and sends $L$ to the graph of the differential of a compactly supported function $S:\R\to\R$. In this case, the Hofer distance is less than $\max S-\min S$. This value can be interpreted as the difference between the areas of some domains delimited by $L$ and $L_0$ respectively in $D_+$ and $D_-$. This leads to the estimate $d(L,L_0)\leq \frac12$, which proves the theorem in the case $\mumax(L)-\mumin(L)=1$.

Suppose now that $\mumax(L)-\mumin(L)>1$. Then, at least one of the two integers $\mumin(L)$ and $\mumax(L)$ does not vanish. Suppose $\mumin(L)< 0$. We are going to construct a Hamiltonian diffeomorphism $h$, with energy $\|h\|\leq 1$, such that $\mumin(h(L))=\mumin(L)+1$. If $\mumax{(L)}>0$, an analogous construction could also be performed to reduce $\mumax{(L)}$ by 1. In any case, we are able to reduce the maximal index gap $\mumax-\mumin$ by 1. It follows by induction that there exists a Hamiltonian diffeomorphism $g$ that sends $L$ to a diameter $L'$ with $\mumax(L')-\mumin(L')=1$ and $\|g\|\leq\mumax(L)-\mumin(L)-1$.
This will conclude the proof of Theorem \ref{theoreme hofer-maslov} since
$$d(L,L_0)\leq d(L,L')+d(L',L_0)\leq \|g\|+\frac12.$$

Let us now suppose $\mumin(L)<0$ and construct a Hamiltonian diffeomorphism $h$ with energy less than 1 and such that $\mumin(h(L))=\mumin(L)+1$.

\begin{figure}
\begin{center}
\includegraphics[width=9cm]{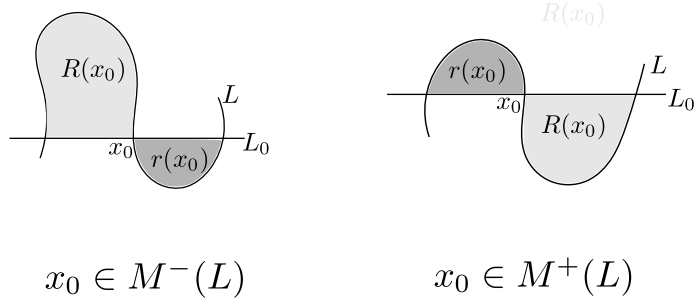} 
\end{center}
\caption{Intersection points belonging to $M^-$ or $M^+$.}
\label{dessin mplus mmoins} 
\end{figure}

%
%
Let $M(L)$ be the (finite) set of transverse intersection points of $L$ with $L_0$ whose index is $\mumin(L)$. Since $\mumin(L)< 0$, the points $x_1$ and $x_N$ do not belong to $M(L)$. 
Let $M^-(L)$ (resp. $M^+(L)$) be the set of all $x\in M(L)$ whose small region is contained in the lower (resp. upper) half-disk (Figure \ref{dessin mplus mmoins}). The points of $M(L)$ are ordered by their position on $L_0$. We will denote by $\preceq$ this order relation and by $\prec$ the associated strict order. Namely, $x\prec y$ means that $x$ is strictly on the left of $y$. In a similar way, a region delimited by $L$, $L_0$ and two intersection points $y,z$ (i.e., a Whitney disk) will be said to be \emph{on the left} of an intersection point $x$ if $y$ and $z$ are on the left of $x$.

When we consider a minimum point in $(M^-(L), \preceq)$, we can refine Lemma \ref{lemme technique} as follows.

\begin{lemma}\label{lemme comb}
Let $x_0$ be the minimum in $(M^-(L),\preceq)$. Then there exists a Hamiltonian diffeomorphism $\phi$ meeting the conditions of Lemma \ref{lemme technique} and satisfying moreover the following:\begin{itemize}
 \item[(i)] $M(\phi(L))\varsubsetneq M(L)$ and $\forall q\in M^-(\phi(L)), q\succ x_0$. 
 \item[(ii)] For any $q\in M(\phi(L))$, consider the region $\mathcal{R}$ of $q$ that is included in $D^-$ where $q$ is seen as a point of $M(L)$. If $\mathcal{R}$ is on the right of $r(x_0)$, then $\mathcal{R}$ does not meet the support of $\phi$.
 \item[(iii)] $M^+(\phi(L))\subset M^+(L)$.
\end{itemize}
\end{lemma}

Let us postpone the proof of this lemma to the end of the section and achieve the proof of Theorem \ref{theoreme hofer-maslov}. 

Applying Lemma \ref{lemme comb} to the minimum point $x_1$ in $(M^-(L),\preceq)$, we get a Hamiltonian diffeomorphism $\phi_1$ meeting several conditions. Then, we repeat inductively the operation on the successive minimum points $x_1, x_2, x_3,\ldots$, constructing diffeomorphisms $\phi_1,\phi_2,\phi_3,\ldots$ until $M^-$ is empty.
The process stops at some point since according to property (i) in Lemma \ref{lemme comb} the (finite) cardinal of $M(L)$ decreases at each step.
We have constructed a finite family of Hamiltonian diffeomorphisms $\phi_1,\ldots,\phi_p$ such that if we write $\psi=\phi_p\circ\ldots\circ\phi_1$ then $M^-(\psi(L))=\emptyset$.

\begin{figure}
\begin{center}
\includegraphics[width=5cm]{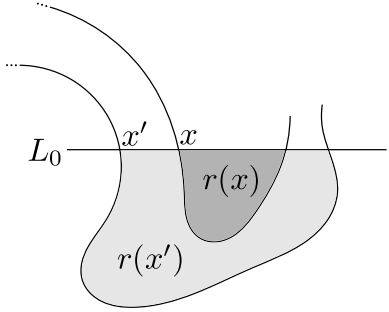} 
\end{center}
\caption{If $r(x)\subset r(x')$, then $x'\preceq x$.}
\label{dessin inclusion petites regions} 
\end{figure}

Moreover, the energy of each $\phi_i$ is approximately the area $w(x_i)$ of the corresponding small region $r(x_i)$. But from conditions (i) and (ii) in Lemma \ref{lemme comb}, the small regions of the points involved in the process are pairwise disjoint. Indeed, these conditions imply that each $r(x_i)$, which is a priori a region of $x_i$ seen as a point in $M(\phi_{i-1}\ldots\phi_1(L))$, is actually already a region of $x_i$ seen as a point in $M(L)$ (the region is invariant under each successive diffeomorphism). But, two such regions are either disjoint or nested. This last possibility
cannot happen for two small regions involded in our process: the largest of those two small regions would be associated to the $\preceq$-smallest intersection point (Figure \ref{dessin inclusion petites regions}), thus we would apply Lemma \ref{lemme comb} first to this larger small region during the process and this would remove the other intersection point. As a consequence, the energy of $\psi$ is bounded above as follows:
$$\|\psi\|\leq\sum_{i=1}^p\|\phi_i\|\leq \sum_{i=1}^p(w(x_i)+\eps)\leq \text{area}(D^-)+p\eps=\frac12+p\eps.$$
Since we can do it for any small enough $\eps>0$, we get $\|\psi\|\leq\frac12$.

Now we apply the same process to the points of $M^+(L)$: we exchange the roles of $M^+(L)$ and $M^-(L)$, which means in particular that we have to start from the right. We get a Hamiltonian diffeomorphism $\psi'$ with Hofer energy less than $\frac12$ and such that $M^+(\psi'\circ\psi(L))=\emptyset$. Moreover, it follows from Lemma \ref{lemme comb} (iii) that $M^-(\psi'\circ\psi(L))\subset M^-(\psi(L))=\emptyset$. Hence $M(\psi'\circ\psi(L))=\emptyset$. 

As a conclusion, we can move with energy at most 1 any diameter with given $\mumin$ to a diameter with minimal index $\mumin+1$. As explained in the beginning of the section, this concludes the proof of Theorem \ref{theoreme hofer-maslov}.
%
%
$\quad\Box$

\begin{demo}[of Lemma \ref{lemme comb}]

\medskip
\noindent(i)
According to Lemma \ref{lemme technique}, $\phi$ removes some intersection points letting the remaining intersection points at the same place and with the same index. Thus, $M(\phi(L))\subset M(L)$ is obvious.

Let $q\in M^-(\phi(L))$. First, $q\neq x_0$ because the intersection $x_0$ is removed by $\phi$. If $q\in M^-(L)$ then $q\succ x_0$ by definition of $x_0$.
Let us show that $\phi$ may be chosen so that $q\in M^+(L)$ implies $q\succ x_0$. The two conditions $q\in M^+(L)$ and $q\in M^-(\phi(L))$ mean that either the area of $r(q)$ is expanded by $\phi$ or that of $R(q)$ is reduced by $\phi$ (maybe both). We will study these two cases separately. But we must first remark that the vector field generating $\phi$ is oriented downward for points of $L_0$ located on the left of (a neighbourhood of) $x_0$ and is oriented upward for points of $L_0$ located on the right of $x_0$ (see Figures \ref{dessin preuve lemme technique} and \ref{dessin hamiltonien simple}).

Suppose now that $r(q)$ is expanded by $\phi$. Then, since it is included in $D^+$, it cannot be entirely on the left of $x_0$ because otherwise, $\phi$ would decrease it. Since  $q$ is on the right of $r(q)$, it follows that $q$ is on the right of $x_0$.

Suppose now that $R(q)$ is reduced by $\phi$. Then a similar argument implies that $R(q)$ is not entirely on the left of $x_0$ and two cases could happen: either $R(q)$ is entirely on the right of $x_0$, or $x_0$ belongs $R(q)$. The first case implies $q\succ x_0$ as wished.  
In the second case, there are again two possiblities (Figure \ref{dessin dernieres possibilites}): (a) $q\notin R(x_0)$, (b) $q\in R(x_0)$. We want to avoid these to cases. To avoid the case (a) we only have to choose the map $\phi$ to have a support sufficiently n so that
$R(q)$ remains invariant. Avoiding the case (b) requires a bit more work and we will need to apply the refined Lemma \ref{lemme technique 2} to appropriate regions $A_i$, $A'_j$ and $B_j$, $1\leq i\leq p$, $0\leq j\leq n$. 

\begin{figure}
\begin{center}
\includegraphics[width=10cm]{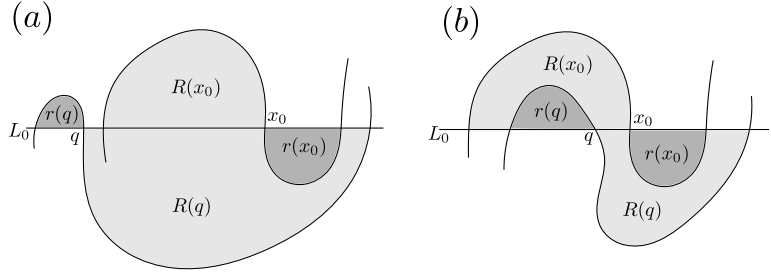} 
\end{center}
\caption{Two possibilities for $x_0\in R(q)$: (a) $q\notin R(x_0)$, (b) $q\in R(x_0)$.}
\label{dessin dernieres possibilites} 
\end{figure} 

Suppose that the set $C$ of points $q\in R(x_0)\cap M^+(L)$ such that $r(x_0)\subset R(q)$ is not empty. We want to choose $\phi$ so that $C$ does not meet $M^-(\phi(L))$.
To do so, we classify the points in $C$ into a partition of four sets $C_1$, $C_2$, $C_3$ and $C_4$, constructed by the following algorithm.

We initialize our sets as $C_1=\emptyset$, $C_2=\emptyset$, $C_3=\emptyset$ and $C_4=C$. We also initialize some auxiliary variables: a set $E=C$, a real variable $\alpha=0$ and a point $q=\max(E)$. Then we iterate the following operations until $E=\emptyset$.
\begin{itemize}
 \item[(1)] If $E\neq \emptyset$ and $\alpha+\omega(q)<w(x_0)$, then modify all the variables by the rules (the symbole $":="$ means "is replaced by"):
\begin{trivlist}
 \item $C_1:= C_1\cup\{q\}$
 \item $C_2:= C_2$
 \item $C_3:=C_3\cup \{q'\in E\,|\,q'\in r(q), q'\neq q\}$
 \item $C_4:=C_4\backslash \{q'\in E\,|\,q'\in r(q)\}$
 \item $\alpha:=\alpha+w(q)$
 \item $E:=E\backslash\{q'\in E\,|\,q'\in r(q)\}$.
\end{trivlist}
 \item[(2)] If $E\neq \emptyset$ and $\alpha+\omega(q)>w(x_0)$, then modify the variables by the rules:
\begin{trivlist}
 \item $C_1:= C_1$
 \item $C_2:= C_2\cup\{q\}$
 \item $C_3:=C_3$
 \item $C_4:=C_4\backslash \{q\}$
 \item $\alpha:=\alpha$
 \item $E:=\{q'\in E\,|\,r(q')\subset r(q), q'\neq q\}$.
\item  
\end{trivlist}
\end{itemize}

The result of the algorithm is illustrated by Figure \ref{dessin algorithm}.
\begin{figure}
\begin{center}
\includegraphics[width=10cm]{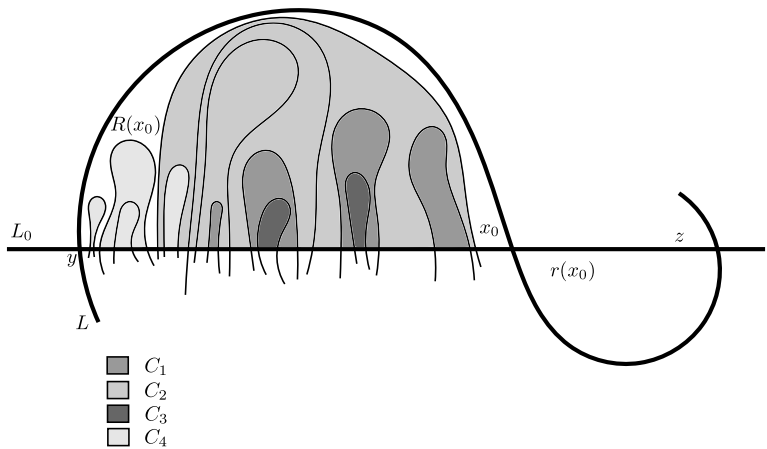} 
\end{center}
\caption{Small regions of points in $C$ are represented. The small regions are colored according to their belonging to $C_1$, $C_2$, $C_3$ or $C_4$ at the end of the algorithm.}
\label{dessin algorithm} 
\end{figure} 
Once it has been performed, we decide to call $A_1,A_2,...$ the small regions of the points of $C_1$, ordered from the right to the left.
We denote by $A'_1, A'_2, \ldots$ the small regions of the points of $C_2$, ordered by inclusion from the smallest to the largest. We also denote by $B_1, B_2,\ldots$ the large regions associated to $A'_1, A'_2, \ldots$. This family of regions satisfies the settings of Lemma \ref{lemme technique 2} which provides us with a diffeomorphism $\phi$ with good properties. Let us analyse how it acts on the set $C$. First, points in $C_1$ are removed. Then, points in $C_2$ remain but the area of their small region decreases by exactly the same amount as their large region. Points in $C_3$ have their small regions included in small regions of points in $C_1$ and thus are removed. Finally points in $C_4$ have their small and large regions disjoint from the support of $\phi$ (they are on the left of its support). In any cases, we see that no points of $C$ belong to $M_-(\phi(L))$, which is exactly what we wished for.

\medskip
\noindent (ii) This condition is obviously satisfied provided the support of $\phi$ is chosen sufficiently narrow.

\medskip
\noindent (iii) 
Let $x\in M^+(\phi(L))$. Suppose $r(x)\subset D^-$ where $x$ is seen as an intersection point of $L$ and $L_0$ (i.e. $x\in M^-(L)$). Then either the region $r(x)$ is increased by $\phi$ or $R(x)$ is decreased (it may be both). But this imposes $x\prec x_0$. Contradiction.
\end{demo}

\end{document}